\documentclass[reqno,12pt]{amsart}

\usepackage[margin=2.5cm]{geometry}
\usepackage[T1]{fontenc}
\usepackage[latin1]{inputenc}
\usepackage[english]{babel}
\usepackage{amsmath,amssymb,amsthm}
\usepackage{enumitem}
\usepackage{mymacros}
\usepackage{hyperref}
\usepackage{xcolor}
\usepackage{graphicx}

\renewcommand{\MR}[1]{}

\title[Holomorphic semigroups and Sarason's characterization of VMOA]{Holomorphic semigroups and Sarason's characterization of vanishing mean oscillation}

\author{Nikolaos Chalmoukis}
\address{Dipartimento di Matematica, Universit\`a di Bologna, 40126, Bologna, Italy}
\email{nikolaos.chalmoukis2@unibo.it}
\thanks{The first author is a member of INdAM}

\author{Vassilis Daskalogiannis}
\address{Department of Mathematics, Aristotle University of Thessaloniki, 54124, Greece }
\email{vdaskalo@math.auth.gr}
\thanks{}


\begin{document}

\begin{abstract}
It is a classical theorem of Sarason that an analytic function of bounded mean oscillation ($BMOA$), is of vanishing mean oscillation if and only if its rotations converge in norm to the original function as the angle of the rotation tends to zero. In a series of two papers Blasco et al. have raised the problem of characterizing all semigroups of holomorphic functions $(\varphi_t)$ that can replace the semigroup of rotations in Sarason's Theorem. We give a complete answer to this question, in terms of a logarithmic vanishing oscillation condition on the infinitesimal generator of the semigroup $(\varphi_t)$. In addition we confirm the conjecture of Blasco et al. that all such semigroups are elliptic. We also investigate the analogous question for the Bloch and the little Bloch space and surprisingly enough we find that the semigroups for which the Bloch version of Sarason's Theorem holds are exactly the same as in the $BMOA$ case.
\end{abstract}

\maketitle

\section{Introduction \& Main results}

A {\it semigroup of analytic self maps} of the unit disc $\bD$ is the flow of a (unique) holomorphic vector field $G$ on $\bD$ which is defined for all positive times  $ t \geq 0 $. In other words it is the solution $(\varphi_t: \bD \mapsto \bD )_{t\geq 0}$ of the Cauchy problem 
\begin{equation} \begin{cases} G(\varphi_t(z))  = \dfrac{\partial \varphi_t(z)}{\partial t};  \\ 
\varphi_0(z)  \equiv z, \end{cases} \end{equation}
when this exists. In particular $G$ is called the {\it infinitesimal generator } of the semigroup. From the dynamical viewpoint the continuous version of Denjoy-Wolff's Theorem \cite[Theorem 8.3.1]{Bracci2020} guarantees the existence of a unique point $\tau \in \overline{\bD}$ which  is called the {\it Denjoy-Wolff point} of $(\varphi_t)$ and $\varphi_t $ converges to $\tau$ uniformly on compact sets as $t \to + \infty$, except when $(\varphi_t)$ consists of elliptic automorphisms of $\bD$ . This allows for a first classification of semigroups; in the {\it elliptic} ones, when $\tau \in \bD$ and the {\it non-elliptic}, when $|\tau| = 1.$ 
 
There has been an increasing amount of literature on semigroups of analytic functions studying both the dynamical features of semigroups, such as slopes of orbits \cite{Betsakos2015, Bracci2020}, rates of convergence to the Denjoy-Wolff point \cite{Betsakos2020} and boundary fixed points \cite{Contreras2004} to name a few, but also studying holomorphic semigroups in the context of (holomorphic) function spaces. This approach was pioneered by Berkson and Porta \cite{Berkson1978}. In this paper we will focus on some questions regarding this latter aspect.

Let $\cH(\bD) $ be the Frech\'et space of holomorphic functions in the unit disc. In \cite{Berkson1978}, Berkson and Porta considered the {\it semigroup of composition operators} 
\[ C_t(f):= f \circ \varphi_t, \quad f\in \cH(\bD) \]
associated to a semigroup $(\varphi_t)$. They proved that each $(\varphi_t)$ induces a {\it strongly continuous semigroup of operators} \cite[Section 34]{lax2002functional} on the classical Hardy space of analytic functions $H^p, p>0 $, i.e. that each operator $C_t$ is a bounded linear operator on $H^p$, the semigroup identity is satisfied, $C_t \circ C_s = C_{t+s}$ and furthermore $C_t$ converges to the identity operator in the strong operator topology as $t\searrow 0.$ Their work has been quite influential, and has naturally lead to analogous considerations in a variety of spaces of analytic functions in the unit disc, among them the classical weighted Bergman spaces \cite{Siskakis1987} and the Dirichlet space \cite{Siskakis1996}. It turns out that the original results of Berkson and Porta continue to hold in these different settings virtually invariable. That is, the composition operators $C_t$ associated to any given semigroup of analytic functions $(\varphi_t)$ form a strongly continuous semigroup of composition operators in all these spaces. 

A different phenomenon arises when one considers some of the most well known {\it non-separable} spaces of analytic functions. The first one to notice this, although not in the language of semigroup theory, was Sarason \cite{Sarason1975}, in the setting of $BMOA$; the space of analytic functions of bounded mean oscillation. For more background on these spaces the reader is referred to Section \ref{sec:semi}.
\begin{thmx}[{\bf Sarason's Theorem}]
Let $\rho_t(z)=e^{it}z, t\geq 0$ be the family of rotations in the unit disc. Then for a function $f\in BMOA$,  \[\lim_{t\searrow 0} \norm{f\circ \rho_t-f}_{BMOA} = 0  \] if and only if $f$ is of vanishing mean oscillation $(VMOA)$.
\end{thmx}
As a matter of fact the rotations $(\rho_t)_{t\geq 0}$, form a semigroup of analytic functions and Sarason's Theorem shows that the composition semigroup induced by  $(\rho_t)$ is {\it not} strongly continuous, rather there exists a {\it maximal closed subspace} of $BMOA$ on which rotations induce a strongly continuous composition  semigroup. We should mention here that Sarason formulates his theorem in the space of functions of bounded mean oscillation in the real line, but as he notices \cite[p.1]{Sarason1975} the result reported here is an equivalent reformulation of his. 

Motivated by this observation, in a series of two papers \cite{Blasco2013, Blasco2008}, Blasco, Contreras, D{\'{\i}}az-Madrigal, Mart{\'{\i}}nez, Papadimitrakis and Siskakis studied composition semigroups in $BMOA$ and in the Bloch space $ \cB $ as long as in their ``little-oh'' versions, $VMOA$ and $\cB_0$. It turns out that strong continuity depends on the characteristics of each specific semigroup $(\varphi_t)$, which led the authors to introduce the {\it maximal subspace of strong continuity}
\[
[\varphi_t,\,X]:= \{f \in X:\;\lim_{t\searrow 0} \norm{C_t(f) - f}_X = 0\}\, ,
\]
i.e. the maximal linear subspace on which $(\varphi_t)$ induces a strongly continuous composition semigroup. It can be proven that when $X=BMOA$ or $\cB$, for all semigroups this is a closed subspace \cite[Proposition 1]{Blasco2008} of $X$. Furthermore,  each $(\varphi_t)$ generates a strongly continuous semigroup $(C_t)$ in $VMOA$ and $\cB_0$, hence in general we have
\begin{equation} \label{eq:inclusions}
X_0 \subseteq [\varphi_t,\,X] \subseteq X,
\end{equation}
where $X_0=VMOA$ or $\cB_0$. \footnote{In the rest of the paper we shall use the shorthand $X$ and $X_0$ to mean that $X$ is the Bloch space or  $BMOA$ and $X_0$ is either the little Bloch space or $VMOA$ respectively.}

In view of Sarason's Theorem, it is quite natural to ask about those $(\varphi_t)$ for which the maximal subspace of strong continuity is {\it minimal}, in the sense of equation \eqref{eq:inclusions}.  In other words we want to find a characterization of semigroups for which $[\varphi_t, X] = X_0$. 

Blasco et al. \cite{Blasco2008} prove that the following ``logarithmic vanishing Bloch'' condition on the infinitesimal generator 
\begin{equation} \tag{LVB} \label{eq:lbb}
    \lim_{|z| \nearrow 1} \frac{1-|z|^2}{G(z)}\log\frac{1}{1-|z|^2} = 0,
\end{equation}
is sufficient so that $[\varphi_t, \cB] = \cB_0$ holds. On the opposite direction they prove that $[\varphi_t, \cB] = \cB_0$ implies \eqref{eq:lbb} under the {\it a priori} hypothesis 
\begin{equation}\label{eq:logbloch} \limsup_{|z| \nearrow 1} \frac{1-|z|^2}{|G(z)|}\log\frac{1}{1-|z|^2} < + \infty. \end{equation}
It should be noted that non-elliptic semigroups always fail to satisfy \eqref{eq:logbloch} therefore the theorem provides no information on the non-elliptic case. Furthermore there exist non trivial elliptic semigroups that also fail to satisfy \eqref{eq:logbloch}. 

In the same work the authors also investigate the case $[\varphi_t, BMOA] = VMOA $. The sufficient condition they obtain is a variation of the logarithmic vanishing Bloch condition adapted to the nature of $BMOA$. For obvious reasons we shall call it ``logarithmic vanishing mean oscillation'' condition; 
\begin{equation} \tag{LVMO} \label{eq:lvmo}
    \lim_{|a| \nearrow 1} \Big(\log\frac{e}{1-|a|^2} \Big)^2 \int_{\bD} \frac{1-|\phi_a(z)|^2}{|G(z)|^2} dm(z) = 0,
\end{equation}
where $\phi_a(z) : = (a-z)/(1-\bar{a}z)$ and $dm$ is the normalized Lebesgue measure on $\bD$. Similarly to the Bloch case the necessity of this condition is proved under the assumption that 
\begin{equation} \label{eq:lbmo}
     \limsup_{|a| \nearrow 1} \Big(\log\frac{e}{1-|a|^2} \Big)^2 \int_{\bD} \frac{1-|\phi_a(z)|^2}{|G(z)|^2} dm(z) < + \infty.
\end{equation}
It is known, and quite straightforward to verify, that $\eqref{eq:lvmo} \implies \eqref{eq:lbb}$, hence the sufficient condition in the $BMOA$ case is apparently stronger than the one for the Bloch space and analogously $ \eqref{eq:lbmo} \implies \eqref{eq:logbloch} $. Hence, a fortiori, all non-elliptic semigroups and some elliptic ones fail to satisfy $ \eqref{eq:lbmo} $.

In view of the above results it is unclear whether there exist non-elliptic semigroups such that the maximal subspace is minimal (either in the Bloch space or in $BMOA$). This problem has been already  posed as a question in \cite[Question 2]{Blasco2013}. 
In this direction the authors in \cite{Blasco2013} provide some necessary conditions for the minimality of the maximal subspace. Suppose that $(\varphi_t)$ is a non-elliptic semigroup with Denjoy-Wolff point $\tau.$ Then Berkson-Porta's formula \cite[Theorem 10.1.10]{Bracci2020} gives the following representation of $G$;
\[ G(z) = (z-\tau)(\bar{\tau}z-1)p(z), \]
where $p$ is a holomorphic function of non-negative real part. Therefore $p$ has a Herglotz representation by some non-negative Borel measure $\mu$, supported on $ \partial \bD.$ In \cite[Corollary 5]{Blasco2013} the authors prove that if $[\varphi_t,X]=X_0$ then $\mu$ has no atoms, i.e. $\mu\{\zeta\}=0, \forall \zeta \in \partial \bD$. Furthermore, they prove \cite[Corollary 6]{Blasco2013}  that if $(\varphi_t)$ is non-elliptic and $[\varphi_t,X]=X_0$ then the Koenigs function $h$ satisfies 
\begin{equation}\label{not BMOA}
h\in \left(\bigcap_{p<\infty}H^p\right)\setminus BMOA,
\end{equation}
where $H^p$ is the classical Hardy space.

We have been able to answer these questions, providing a complete characterization of the semigroups for which $[\varphi_t, X]=X_0$.
\begin{thm}\label{main:thm:elliptic}
Let $ (\varphi_t) $ be a semigroup of analytic functions with infinitesimal generator $G$. The following are equivalent; 
\begin{itemize}
    \item[(a)] $ [\varphi_t, \cB]=\cB_0$,
  \vspace{0.1cm}
    \item[(b)]  $ [\varphi_t, BMOA]=VMOA $,
   \vspace{0.1cm}
    \item[(c)] $(\varphi_t)$ is an elliptic semigroup and $G$ satisfies the logarithmic vanishing Bloch condition \eqref{eq:lbb},
      \vspace{0.1cm}
    \item[(d)] $(\varphi_t)$ is an elliptic semigroup and $G$ satisfies the logarithmic vanishing mean oscillation condition \eqref{eq:lvmo}. 
\end{itemize}
\end{thm}

The surprising aspect of this theorem is that not only the sufficient conditions of Blasco et al. are also necessary for the minimality of the maximal subspace under no further assumptions, but quite unexpectedly the two conditions are equivalent. Hence the class of holomorphic semigroups which can replace the rotations in Sarason's theorem is exactly the same for the Bloch space and for $BMOA$.  In particular the implications $(a) \implies (c)$ and $(b) \implies (d)$  answer in the affirmative the question of Blasco et al. \cite{Blasco2013}, i.e. that no non-elliptic semigroup has maximal subspace which coincides with either $\cB_0$ or $VMOA$.

\subsection*{Plan of the paper} In Section \ref{sec:semi} we give a quick overview of the concepts that will go into the proof of the main theorem. In particular we discuss in more detail semigroups of analytic functions and the spaces $BMOA$ and $\cB$ as well as some weighted versions of them. In Section \ref{sec:mainproofs} we prove the main theorem. In fact, the central part of the proof is a construction presented in Proposition \ref{prop:mainconstruction}.

\subsection*{Acknowledgments} The authors are indebted to professors Aristomenis Siskakis and Petros Galanopoulos for very helpful discussions. We would like also to thank the anonymous referees for carefully reading the manuscript and providing helpful comments.

\section{Background}\label{sec:semi}
In this section we shall discuss some of the background material and introduce some notation that we are going to use later. In the unit disc $\bD$ we denote by $\delta$ the hyperbolic distance,
\[
\delta(a,z):=\dfrac{1}{2}\log\dfrac{1+|\phi_a(z)|}{1-|\phi_a(z)|},
\quad \text{where,} \quad \phi_a(z):=\frac{a-z}{1-\bar{a}z} \;\;a,z\in \bD.
\]
This is the distance corresponding to the hyperbolic Riemannian metric $ds/(1-s^2)$. The metric space $(\bD, \delta)$ is a model of the hyperbolic plane usually called the {\it Poincar\'e disc}. The functions $\phi_a$ are   isometric automorphisms of the Poincar\'e disc and are also involutions ($\phi_a^{-1} = \phi_a $). 

For a holomorphic function $f$ defined on $\bD$ we define its {\it hyperbolic translation} $f_a$ with respect to $a\in \bD$ as 
\[ f_a(z):= f(\phi_a(z))-f(a). \]
Another fact that is going to be used repeatedly is the following approximation for the hyperbolic distance of a point $z\in \bD$ to the origin; 
\[ 1+\delta(0,z) \approx \log\frac{e}{1-|z|^2}. \]


Let us now take a closer look to holomorphic semigroups. An equivalent way to define a holomorphic semigroup $(\varphi_t)$ is as a family $\{\varphi_t:\,t\geq0\}$ of analytic self maps of the unit disc $\varphi_t : \bD \to \bD$ such that
\begin{enumerate}
    \item $\varphi_0(z)\equiv z$
    \item $\varphi_t \circ \varphi_s = \varphi_{t+s},\;\;t,s\geq 0$
    \item $\varphi_t(z) \to z$ uniformly on compact subsets of $\bD$, as $t\searrow 0$.
\end{enumerate}
It turns out that if $(\varphi_t)$ is a semigroup then each $\varphi_t$ is univalent \cite[Theorem 8.1.17]{Bracci2020}. 
In addition, for all members of a semigroup $(\varphi_t)$ (other than the hyperbolic rotations) there exists a common ``fixed point'' $\tau \in \ol{\bD}$, for which
\[
\lim_{t\to\infty}\varphi_t(z)=\tau, \;z\in\bD\,,
\]
usually called the Denjoy-Wolff point of $(\varphi_t)$. The concept of the Denjoy-Wolff point of a semigroup holds a key role in the semigroup theory and we can classify semigroups with respect to their Denjoy-Wolff point, $\tau$, as follows \cite[Theorem 8.3.1]{Bracci2020};
\begin{enumerate}
\item If $\tau\in \bD$, then there exists $\lambda\in \bC\setminus \{0\}$ with $\Re(\lambda)\geq 0$ such that
\[
\varphi_t'(\tau)=e^{-\lambda t},\;\;t\geq 0\,.
\]
In addition, we either have $|\varphi_t'(\tau)|=1$ for every $t>0$, or $|\varphi_t'(\tau)|<1$ for every $t>0$.
\item If $\tau \in \partial\bD$, then there exists $\lambda \geq 0$ such that
\[
\angle\lim_{z\to \tau}\varphi_t'(\tau)=e^{-\lambda t},\;\;t\geq 0.
\]

\end{enumerate}
The number $\lambda$ is called the spectral value of the semigroup.
We say that $(\varphi_t)$ is {\it elliptic} if $\tau\in\bD$, {\it parabolic} if $\tau\in \partial \bD$ with spectral value $\lambda=0$, and {\it hyperbolic} if $\tau\in \partial \bD$ with spectral value $\lambda>0$. The semigroup is called {\it non-elliptic} if it is either parabolic or hyperbolic.

 If $\,(\varphi_t)\,$ is a semigroup then the limit
\[
G(z)=\lim_{t\searrow 0}\dfrac{\varphi_t(z)-z}{t}
\]
exists uniformly on compact subsets of $\bD$. The function $G\in \cH(\bD)$ is 
the infinitesimal generator of $(\varphi_t)\,$ and characterizes the semigroup in a unique way. In addition, $G$ satisfies the following relations
\begin{equation}\label{generator}
G(\varphi_t(z))\,=\,\dfrac{\partial \varphi_t(z)}{\partial t}\,=\,G(z)\,\dfrac{\partial \varphi_t(z)}{\partial z},\;\;z\in\bD\,,\;\;t\geq 0\,.
\end{equation}
Due to the Berkson-Porta formula \cite{Berkson1978} we can represent the infinitesimal generator $G$ in terms of the Denjoy-Wolff point $\tau$ of the semigroup as
\begin{equation}\label{D-W}
G(z)\,=\,(\bar{\tau}z -1)(z-\tau)p(z),\;\;z\in\bD\,,
\end{equation}
where $\,\tau\in \overline{\bD}\,$ and $\,p\in \cH(\bD)\,$ with $\Re\left(p(z)\right)\geq 0\,$ for all $z\in\bD$. Conversely every function of this form is the infinitesimal generator of a holomorphic semigroup.

A geometric description of all holomorphic semigroups is provided by the so called {\it Koenigs function}, a conformal map which conjugates a  given semigroup $(\varphi_t)$ to a model semigroup. 

When $(\varphi_t)$ is an elliptic semigroup, with Denjoy--Wolff point $\tau\in \bD$ the  function $h$ is the unique conformal map such that $h(\tau)=0,\;h'(\tau)=1$ and 
\[
h(\phi_t(z))=e^{-\lambda t}h(z),\;\;z\in\bD,\;t\geq0\,.
\]
In addition we have that $\dfrac{h'(z)}{h(z)}=-\dfrac{\lambda}{G(z)}$. In the non-elliptic case, $h$ is the unique conformal map such that  $h(0)=0$ and
\[
h(\phi_t(z))=h(z)+it,\;\;z\in\bD,\;t\geq0\,.
\]
In this case we have that $h'(z)=\dfrac{i}{G(z)}$.

For a semigroup $(\varphi_t)$ with infinitesimal generator $G$ and Denjoy - Wolff point $\tau$, following the notation used in \cite[Definition 4]{Blasco2013}, we consider the function $\gamma: \bD \to \bC$, which we will call the {\it associated $\gamma$-symbol} of $(\varphi_t)$. This function is defined as follows. If $\tau\in\bD$, then
\[
\gamma (z):= \int_\tau^z \dfrac{\zeta -\tau}{G(\zeta)}\,d\zeta \,,
\]
while if $\tau\in \partial \bD$, then
\[
\gamma(z):=\int_0^z \dfrac{i}{G(\zeta)}\,d\zeta \,.
\]
 In the case where $\tau\in \partial \bD$, then $\gamma$ coincides with $h$, while if $\tau\in\bD$, then $\gamma'(z)\,=\, -\frac{z-\tau}{\lambda}\frac{h'(z)}{h(z)}$.

The maximal subspace of strong continuity, for a semigroup of composition operators $(C_t)$, can also be described in terms of the infinitesimal generator $G$ \cite[Theorem 1]{Blasco2013}. If $(C_t)$ acts on a Banach space $\cX$ of analytic functions in the unit disc which contains the constant functions, and in addition we have that $\sup\limits_{t\leq 1} \norm{C_t}_\cX\,<\,\infty$, then
\begin{equation}\label{MaxSubG}
    [\varphi_t,\,\cX]\,=\,\ol{\{f \in \cX:\;Gf' \in \cX\}}.
\end{equation}

This description already indicates a connection between the maximal subspace and the so called {\it generalized Volterra operator} defined for an analytic symbol $g$ as
\[
T_g(f)(z)\,:=\,\int_0^z f(\zeta)g'(\zeta)\,d\zeta\,,\;f \in \cH(\bD) \,.
\]

This operator was first introduced by Pommerenke \cite{Pommerenke1977}, who studied its boundedness properties on the Hardy space $H^2$ in connection to the analytic John-Nirenberg inequality. Since then, several authors studied these operators focusing on conditions on the symbol $g$ under which $T_g$ is bounded or compact. The survey papers \cite{ Aleman2006, Siskakis2004} contain much more information on the generalized Volterra operator.

It turns out that if $\cX$ is Banach space of analytic functions,  and  $g$ is the associated $\gamma$-symbol of $(\varphi_t)$, under mild additional assumptions on $\cX$, we have the following characterization for the maximal subspace of strong continuity \cite[Proposition 2]{Blasco2013};
\begin{equation}\label{MaxSubTg}
[\varphi_t,\,\cX]=\overline{\cX\cap (T_\gamma(\cX)\oplus \cC)}\,,
\end{equation}
where $\cC$ is the set of all constant functions.

Finally we introduce some definitions and we recall some theorems regarding the  Banach spaces of analytic functions we are interested in. The space $BMOA$ is the space of all analytic functions in the Hardy space $H^2$, which have bounded mean oscillation. Having to choose between many equivalent descriptions, we will use a description in terms of Carleson measures. We say that $f\in BMOA$, if and only if
\[
\norm{f}^2_*:=\sup_{I\subseteq \partial\bD }\, \frac{1}{|I|}\,\int_{S(I)} \vert f'(z)\vert^2(1-|z|^2)\,dm(z)\,<\,\infty\,,
\]
where $dm(z)=\frac{dx\,dy}{\pi}$ is the normalized area Lebesgue measure in $\bD$, $I$ is any arc on $\partial\bD$ and $\vert I \vert$ is its length. Also, $S(I)$ is the so called {\it Carleson box}, which for us will be the closed hyperbolic halfplane in the Poincar\'e disc which has $I$ as its boundary. There exist more ``square'' versions of Carleson boxes but the invariant nature of this definition will simplify some of our computations. 

The space $BMOA$ is a Banach space, equipped with the norm
\[
\norm{f}_{BMOA}\,:=\,\vert f(0) \vert \,+\, \norm{f}_*\,.
\]
The closure of all polynomials in $BMOA$ is the space $VMOA$ which has an equivalent description in terms of the following vanishing Carleson condition; 
\[
\lim_{\vert I \vert \searrow 0}\,\frac{1}{|I|}\,\int_{S(I)} \vert f'(z)\vert^2(1-|z|^2)\,dm(z)\,=\,0\,.
\]

The space $BMOA$ is a subspace of the well known Bloch space, denoted by $\cB$. We say that a function $f\in \cH(\bD)$ belongs to $\cB$, if and only if
\[
\sup_{z \in \bD} \vert f'(z) \vert\,(1-\vert z \vert^2)\,<\,\infty\,.
\]
The closure of polynomials in the Bloch norm, is called the little Bloch space, denoted by $\cB_0$. Equivalently  $f\in \cB_0$ if and only if
\[
\lim_{\vert z \vert \to 1} \vert f'(z) \vert\,(1-\vert z \vert^2)=0\,.
\]
The spaces $\cB$ and $\cB_0$, are Banach spaces equipped with the norm
\[
\norm{f}_{\cB}\,:=\,\vert f(0) \vert \,+\, \sup_{z \in \bD} \vert f'(z) \vert\,(1-\vert z \vert^2)\,.
\]

For more information on these spaces see \cite{zhu2007operator}.

In \cite[Definition 3]{Blasco2013}, the authors consider some weighted versions of $BMOA$ and $\cB$ which are closely related to the conditions \eqref{eq:lvmo} and \eqref{eq:lbb} . Let $f\in \cH(\bD)$. Then we say that $f$ belongs to $BMOA_{\log}$ if and only if
    \begin{equation}\label{BMOA_log}
        \sup_{I\subseteq \partial\bD }\, \frac{\left(\log\dfrac{e}{|I|}\right)^2}{|I|}\,\int_{S(I)} \vert f'(z)\vert^2(1-|z|^2)\,dm(z)\,<\,\infty\,,
    \end{equation}
    and $f\in VMOA_{\log}$ if and only if
    \begin{equation}\label{VMOA_log}
      \lim_{\vert I \vert \to 0}\,  \frac{\left(\log\dfrac{e}{|I|}\right)^2}{|I|}\,\int_{S(I)} \vert f'(z)\vert^2(1-|z|^2)\,dm(z)\,=\,0\,.
    \end{equation}
Respectively, we say that $f$ belongs to $\cB_{\log}$ if and only if
    \begin{equation}\label{B_log}
      \sup_{z \in \bD}\, \vert f'(z) \vert\,(1-\vert z \vert^2)\left(\log\dfrac{e}{1-|z|^2}\right)\,<\,\infty\,,  
    \end{equation}
    and $f\in \cB_{\log,\,0}$ if and only if
    \begin{equation}\label{B_log,0}
      \lim_{\vert z \vert \to 1}\, \vert f'(z) \vert\,\left(\log\dfrac{e}{1-|z|^2}\right)\,\,(1-\vert z \vert^2)=0\,. \end{equation}
These spaces naturally appeared in the study of multipliers for $BMOA$ and the Bloch space. A function $g$ is a pointwise multiplier of $\cB$, i.e. $gf\in\cB,\;\;\forall  f\in\cB$ if and only if $g\in H^\infty \cap B_{\log}$ \cite{Brown1991}, where $H^\infty$ is the space of bounded analytic functions in the unit disc. An analogous result holds for $BMOA$, that is; $g$ is a multiplier for $BMOA$, if and only if $g\in H^\infty \cap BMOA_{\log}$ \cite{Ortega1996}.  
  For our purposes these spaces are interesting because they characterize the boundedness and compactness of $T_g$ on $BMOA$ and $\cB$  \cite[Theorems 5 \& 6]{Blasco2013},  \cite{ Siskakis1999, Laitila2011}.
\begin{thmx}\label{bound.Comp.}
$T_g: BMOA\to BMOA$ is bounded, if and only if $\,g\in BMOA_{\log}$. Furthermore, the following are equivalent:
\begin{itemize}
    \item[(i)] $T_g: BMOA\to BMOA$ is compact,
    
    \item[(ii)] $g\in VMOA_{\log}$,
    
    \item[(iii)] $T_g: BMOA\to BMOA$ is weakly compact.
\end{itemize}

An analogous result holds for the Bloch space.
$T_g: \cB\to \cB$ is bounded, if and only if $\,g\in B_{\log}$. Furthermore, the following are equivalent:
\begin{itemize}
    \item[(I)] $T_g: \cB\to \cB$ is compact,
    
    \item[(II)] $g\in \cB_{\log,0}$,
    
    \item[(III)] $T_g: \cB\to \cB$ is weakly compact.
\end{itemize}
\end{thmx}
\noindent In addition, Gantmacher's theorem  \cite[Theorem 5.23]{Aliprantis2006} implies that $T_g: X\to X$ is weakly compact if and only if $T_g(X)\subseteq X_0$.

\subsection*{Notation} For two quantities $A, B$ depending on a number of parameters we shall write $A \lesssim B $ if there exists some positive constant $C>0$ not depending on the parameters such that $A \leq C B.$ The set of parameters in question should be clear from the context. Similarly we shall write $A \approx B$ if $A \lesssim B$ and $ B \lesssim A $.

\section{Proof of main results}\label{sec:mainproofs}

We can now turn towards the proof of Theorem \ref{main:thm:elliptic} which can be divided on a macroscopic scale in two main parts. The first one regards the equivalence of parts $(c)$ and $(d)$ of Theorem \ref{main:thm:elliptic}. This equivalence is in reality another manifestation of the rigidity properties of univalent functions. 

The idea that we employ already appears in \cite{Pommerenke1978} and further refined in \cite{Aulaskari1997}. Here we shall adapt it in the weighted setting relevant to our problem. Although we could do all calculations for the logarithmic weight we prefer to work with a more general class of weights since we think that this renders more clear the idea of the proof. Let $\omega $ be a strictly positive  weight of the class $C^1(\bD).$ We assume the following regularity condition on $\omega$; 
\begin{equation}\label{eq:regularity}
    (1-|z|^2)|\nabla \omega (z) | \leq C_\omega \omega(z), \quad \forall z \in \bD 
\end{equation}for some $C_\omega > 0.$

\begin{lem}\label{lemma:estimateBloch}
Suppose that $\omega$ is a weight which satisfies \eqref{eq:regularity} with some constant $C_\omega < 1. $ Let also $f\in \cH(\bD)$ such that 
\[ |f'(z)| (1-|z|^2) \omega(z) \leq K, \quad \forall z \in \bD, \]
where $K>0$. Then,
\[\int_0^1 \sup_{a \in \bD, |z| \leq r } (\omega(a) |f_a(z)|)^2 dr < + \infty.  \]
\end{lem}

\begin{proof} Let us start with a local oscillation estimate on $\omega$.  Let $z=r e^{i\theta} $
\begin{equation*}
    \log \frac{\omega(z)}{\omega(0)}  \leq \int_0^{r} \frac{|\nabla \omega (se^{i\theta}) | }{ \omega (se^{i\theta}) }ds \leq C_\omega \delta(0,z).
\end{equation*}
Since condition \eqref{eq:regularity} is invariant under composition with M\"obius transformations \cite[Proposition 3.1]{Aleman2009} we have that
\[  \omega(z) \leq e^{C_\omega \delta(z,w)} \omega (w), \quad z,w \in \bD.  \]
We proceed now to an estimate of the quantity that appears in the lemma.
For $z\in \bD$, we have 
\begin{align*}
    |f_a'(z)|(1-|z|^2)  = |f'(\phi_a(z))|(1-|\phi_a(z)|^2)\ 
    \leq K \omega(\phi_a(z))^{-1}.
\end{align*}
    
Hence for $z=re^{i\theta}$ we have,
 \begin{align*}
      \omega(a) |f_a(z)| & \leq K \int_0^r \frac{\omega(a)}{\omega(\phi_a(te^{i\theta}))} \frac{dt}{1-t^2}\leq K \sup_{\{w\, :\, \delta(a,w)\leq \delta(0,r)\} } \frac{\omega(a)}{\omega(w)} \delta(0,r) \\ 
      & \leq K e^{C_\omega \delta(0,r)} \delta(0,r)
 \end{align*}

The lemma follows from the fact that the function 
\[ e^{2 C_\omega \delta (0, r)} \delta(0,r)^2 = \Big( \frac{1+r}{1-r} \Big)^{C_\omega} \Big( \frac12 \log\frac{1+r}{1-r} \Big)^2, \] is integrable in $(0,1) $ if $C_\omega < 1 $.
\end{proof}

The next proposition is a weighted version of Pommerenke's result \cite{Pommerenke1978}.

\begin{prop}\label{prop:univalent}
Let $f:\bD \to \bC$ be univalent and $\omega$ a weight as in Lemma \ref{lemma:estimateBloch}. Suppose also that 
\begin{equation} \label{omega-vanishing} \lim_{|z|\nearrow 1} |f'(z)|(1-|z|^2) \omega (z) = 0 \end{equation} Then,
\begin{equation} \label{omega-bmoa-vanishing} \lim_{|a|\nearrow 1} \omega(a)^2 \int_{\bD}|f'(z)|^2(1-|\phi'_a(z)|^2)dm(z) = 0.\end{equation}
\end{prop}
\begin{proof}
Let $f$ be such a function. Setting $\bD_r=\{z\in\bD: |z|\leq r\},\;r\in[0,1)$, we have
\[
\begin{split}
 \omega(a)^2 & \int_\bD  |f'(w)|^2 (1-|\phi_a(w)|^2)\,dm(w)  
 = \omega(a)^2\int_0^1 \int_{\bD_r} |f'_a(z)|^2 \,dm(z)\,dr \\
 & = \omega(a)^2 \int_0^R \int_{\bD_r} |f'_a(z)|^2 \,dm(z)\,dr 
 +  \omega(a)^2 \int_R^1 \int_{\bD_r} |f'_a(z)|^2 \,dm(z)\,dr 
\\
&=: \, \text{I} \,+\, \text{II}\,.
\end{split}
\]
Since $f_a$ is univalent, the inner integral is the normalized area of the image $f_a(\bD_r)$, hence
\[
\int_{\bD_r} |f'_a(z)|^2 \,dm(z)\leq \, \sup_{z\in \bD_r}|f_a(z)|^2 \,.
\]
Let now $\varepsilon > 0$. By Lemma \ref{lemma:estimateBloch} there exists some $R_0<1$ such that 
\[ \text{II} \leq  \omega(a)^2 \int_{R_0}^1  \sup_{z \in \bD_r } |f_a(z)|^2  dr < \varepsilon, \quad \forall a\in\bD.  \]


In order to estimate integral $\text{I}$, notice first that the oscillation estimate in the proof of Lemma \ref{lemma:estimateBloch}, for $z\leq R_0$, gives
\[
\frac{\omega(a)}{\omega(\phi_a(z))} \leq e^{C_\omega \delta(0,R_0)}, a\in \bD.
\]
Hence we have, 
\[
\begin{split}
\text{I}= & \omega(a)^2 \int_0^{R_0} \int_{\bD_r}\left( \dfrac{|f'(\phi_a(z))|(1-|\phi_a(z)|^2)}{1-|z|^2} \right)^2 \,dm(z)\,dr\\
= &   \int_0^{R_0} \int_{\bD_r} \Big( \frac{\omega(a)}{\omega(\phi_a(z))} \Big)^2 \left( \dfrac{|f'(\phi_a(z))|(1-|\phi_a(z)|^2)\omega(\phi_a(z))}{1-|z|^2} \right)^2 \,dm(z)\,dr \\ 
\leq &  e^{2C_\omega \delta(0,R_0)} \,\int_0^{R_0} 
 \dfrac{1}{(1-r^2)^2} \int_{\bD_r} |f'(\phi_a(z))|^2(1-|\phi_a(z)|^2)^2\omega(\phi_a(z))^2\,dm(z) \,dr.
\end{split}
\]

Now, since $f$ satisfies \eqref{omega-vanishing}, we can find $R_1<1$ such that  
\[ \omega(w)(1-|w|^2)|f'(w)| \leq e^{-C_\omega \delta(0,R_0)}(1-R_0^2)\sqrt{\varepsilon}, \quad \forall w : R_1<|w|<1. \]

Finally there exists some $\delta>0$ such that $|\phi_a(z)|> R_1$, if $|z|\leq R_0$ and $|a| > 1-\delta. $ Which gives
\[
\text{I}\leq  (1-R_0^2)^2 \int_0^{R_0} \dfrac{1}{(1-r^2)^2}\,dr\; \varepsilon\; \leq\;  \varepsilon \,.
\]
Therefore we have proved that for each $\varepsilon > 0$, we can find $\delta>0$ such that for $|a| > 1-\delta$  
\[
\omega(a)^2 \int_\bD  |f'(w)|^2 (1-|\phi_a(w)|^2)\,dm(w) \,\leq\, 2\varepsilon.
\]
\end{proof}

\begin{cor}\label{corollary}
If $f$ is univalent, then 
\[ f\in \cB_{\log,0}\;\; \text{if and only if}\;\; f \in VMOA_{\log}. \]
\end{cor}
\begin{proof}
It is sufficient to prove the direct implication. Consider the weight $ \omega_K(z):=\log\frac{K}{1-|z|^2} $. For some $K>0$ large enough $\omega_K $ satisfies the hypothesis of Lemma \ref{lemma:estimateBloch}. Then $f\in \cB_{\log,0}$ is equivalent to \eqref{omega-vanishing}, hence it satisfies \eqref{omega-bmoa-vanishing} which is equivalent to $f\in VMOA_{\log}$.
\end{proof}
We now turn to the second part of the proof. Roughly speaking the characterization of the maximal subspace of strong continuity by Blasco et al. \cite{Blasco2013} as $\overline{(T_\gamma(X)\oplus \cC) \cap X}$, allows us to approach the problem of studying the maximal subspace of strong continuity purely by functional analytic methods. 
The central part of the proof will therefore follow from a construction of a function in the range of $T_g$ under some technical assumptions on $g$. The construction turns out to be quite explicit by pasting together some holomorphic ``building blocks''. These so called building blocks behave much like the logarithmic function 
 $\ell_w(z):=\log(\frac{e}{1-\bar{w}z})$
in the sense that at a prescribed point (in this case $w$) achieves the biggest possible growth while keeping the $BMOA$ or Bloch norm below a fixed threshold.  Our construction requires some improved decay properties away from the point $w$.  The exact definition of these functions is presented in the next lemma. 

\begin{lem}\label{lem:buildingblock}
For a point $w \in \bD$ we denote by $w^*$ the hyperbolic midpoint between $0$ and $w$. Let also $I_w$ be the closed arc in the unit circle such that $w$ is the point in $S(I_w)$ closest to the origin (see figure \ref{fig:midpoint}).
Then the function 
\[ \beta_w(z) : = \log\frac{e}{1-\phi_{w^*}(z) \overline{w}} \]
satisfies the following properties 
\begin{enumerate}
    \item[(i)] $\norm{\beta_w}_\cB \lesssim \norm{\beta_w}_* \lesssim 1 $,
    \item[(ii)] $\Re \beta_w \geq 0 $,
    \item[(iii)] $|\Im \beta_w | \leq \frac{\pi}{2}$,
    \item[(iv)] $\Re \beta_w(z) \approx \log\frac{e}{1-|w|^2}$ for $z\in S(I_w)$,
    \item[(v)] If $z \not\in S(I_{w^*})$ then, $|\beta_w(z)| \leq c_0$, where $c_0$ is an absolute constant. In particular, for all $\delta>0$ there exists $\delta'>0$ such that if $1-|w|\leq \delta' $ and $1-|z|\geq \delta $, then $|\beta_w(z)| \leq c_0. $
\end{enumerate}
\end{lem}

\begin{proof}
To prove part (i) notice that the M\"obius invariant part of the norm does not change after composing the function $\log\frac{e}{1-\overline{w}z}$ with the M\"obius transformation $\phi_{w^*}$. Also 
\[ \beta_w(0)=\log \frac{e}{1+w^*\overline{w}} \leq 1.
\]
Parts $(ii)$ and $(iii)$ follow by the definition of the logarithm. Then to see $(iv)$ notice that $\phi_{w^*}$ preserves the diameter passing through $w$, it maps $w^*$ at $0$ and it leaves invariant the hyperbolic distances, therefore it should map $S(I_w)$ to $S(I_{w^*}).$ Hence, if $ z\in S(I_w), \, y := \phi_{w^*}(z) \in S(I_{w^*}) $ 
\[ \Re \beta_w(z) = \log\frac{e}{|1-y\overline{w}|} \gtrsim\log\frac{e}{1-|w^*|^2} \approx \log\frac{e}{1-|w|^2}.  \]

\begin{figure}
    \centering
    \includegraphics[scale=0.4]{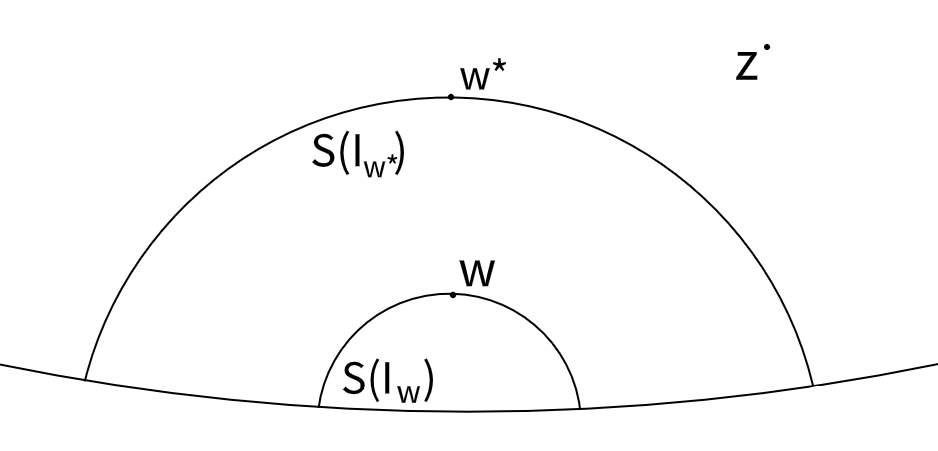}
    \caption{The construction in Lemma \ref{lem:buildingblock}}
    \label{fig:midpoint}
\end{figure}

A similar geometric reasoning as before shows that $\phi_{w^*}(\bD \setminus S(I_{w^*})) $ is the half disc which contains $-w$ and  defined by the diameter perpendicular to the radius passing from $w$. Let therefore $z\in\bD$ in this half plane, or equivalently $\Re(z\overline{w})\leq 0.$

Hence, 
\[ \Big| \log\frac{e}{1-z\overline{w}} \Big| \leq \log\frac{e}{|1-z\overline{w}|} + \frac{\pi}{2} \leq 1+\frac{\pi}{2} \leq 3. \]

To verify the second part of property $(v)$  it remains only to notice that for every $\delta>0 $ there exists $\delta'>0$ such that if $1-|w|\leq \delta'$, then the disc $\{ |z| \leq 1-\delta  \}$ is contained in $\bD \setminus S(I_{w^*}).$
\end{proof}

The next proposition is the main technical tool in the proof of our main theorem.

\begin{prop}\label{prop:mainconstruction}
Let $g\in BMOA \setminus VMOA_{\log}$. Then there exists a function $F\in BMOA$ such that $T_gF\in BMOA \setminus VMOA$.
\end{prop}
\begin{proof} First note that if $g\in BMOA \setminus VMOA $ then the function $F \equiv 1 $ satisfies the required properties. On the other hand 
if $g \in BMOA_{\log} \setminus VMOA_{\log}$ then by Theorem \ref{bound.Comp.} we have that $T_g(BMOA) \subset BMOA$ but $T_g(BMOA) \not \subseteq VMOA$ therefore we can find  a function $F$ as claimed in the thesis of the theorem. 
We have therefore reduced the problem to  the case $g\in VMOA\setminus BMOA_{\log}$.

In order to reduce the number of constants in the proof we assume without loss of generality  that
\[ \int_\bD |g'(z)|^2(1-|z|^2)dm(z) = 1 \leq \norm{g}_*.\]

\textbf{Basic reduction of the problem.} The basic step of the proof is a construction of a sequence of arcs $\{I_n\} $ and a sequence of functions $F_n$ of the form
\[ F_n(z)=\sum_{k=0}^n a_k \beta_{w_k}(z), \]
for some $w_k\in \bD$ such that; 
\begin{enumerate}
    \item The coefficients $\{a_k\}$ satisfy $0\leq  a_k \leq 2^{-k}$.
    \item For all $n\in \bN$ it holds 
    \[ \frac{1}{|I_n|}\int_{S(I_n)}(\Re F_n(z))^2|g'(z)|^2 (1-|z|^2) dm(z) \geq 1. \]
    \item For all $ n \in \bN $ we have $\norm{T_gF_n}_* \leq \max\{\norm{T_g F_{n-1}}_*+2^{-n}C(g), C(g) \}$ where $C(g)$ is a positive constant which depends only on $g$.
\end{enumerate}
Suppose now that we can construct such a sequence of functions. Then we can finish the proof as follows; we have that \[ \sum_{k=0}^\infty a_k( \norm{\beta_{w_k}}_*+|\beta_{w_k}(0)|) \lesssim \sum_{k=0}^\infty a_k < +\infty.  \]
Therefore 
\[ F:= \sum_{k=0}^\infty a_k \beta_{w_k} \in BMOA. \]

Applying repeatedly property $(3)$ we find that 
\[ \norm{T_gF_n}_* \leq \max\{ \norm{T_gF_0}_* + \sum_{r=1}^n 2^{-r} C(g), C(g) \} \leq \norm{g}_*+C(g). \]
There is a slight subtlety in the fact that since $T_g$ is not continuous we cannot directly infer from $(3)$ that $T_gF\in BMOA$. But this problem is easily overcome. It suffices to prove that $ |\Im F(z)g'(z)|^2(1-|z^2|)dm(z) $ and $ (\Re F(z))^2|g'(z)|^2(1-|z|^2)dm(z)$, are Carleson measures for the Hardy space. The first one is clearly a Carleson measure since the imaginary part of $F$ is bounded. Let $I\subseteq \partial \bD$. By the monotone convergence theorem we have 
\begin{align*} \frac{1}{|I|}\int_{S(I)}(\Re F(z))^2|g'(z)|^2(1-|z|^2)dm(z) & = \lim_{n\to \infty}  \frac{1}{|I|}\int_{S(I)}(\Re F_n(z))^2|g'(z)|^2(1-|z|^2)dm(z) \\ & \leq \lim_{n\to \infty }\norm{T_g F_n}_*^2 \leq C(g)^2. 
\end{align*}

Finally we should prove that $T_gF\not\in VMOA.$ This is a simple consequence of (2); 
\[ \frac{1}{|I_n|}\int_{S(I_n)}|F(z)g'(z)|^2(1-|z|^2)dm(z) \geq \frac{1}{|I_n|}\int_{S(I_n)}(\Re F_n(z))^2|g'(z)|^2(1-|z|^2)dm(z) \geq 1. \]

Therefore it remains only to construct such a sequence of functions. This will be done in a recursive way.
\begin{figure}
    \centering
    \includegraphics[scale=0.3]{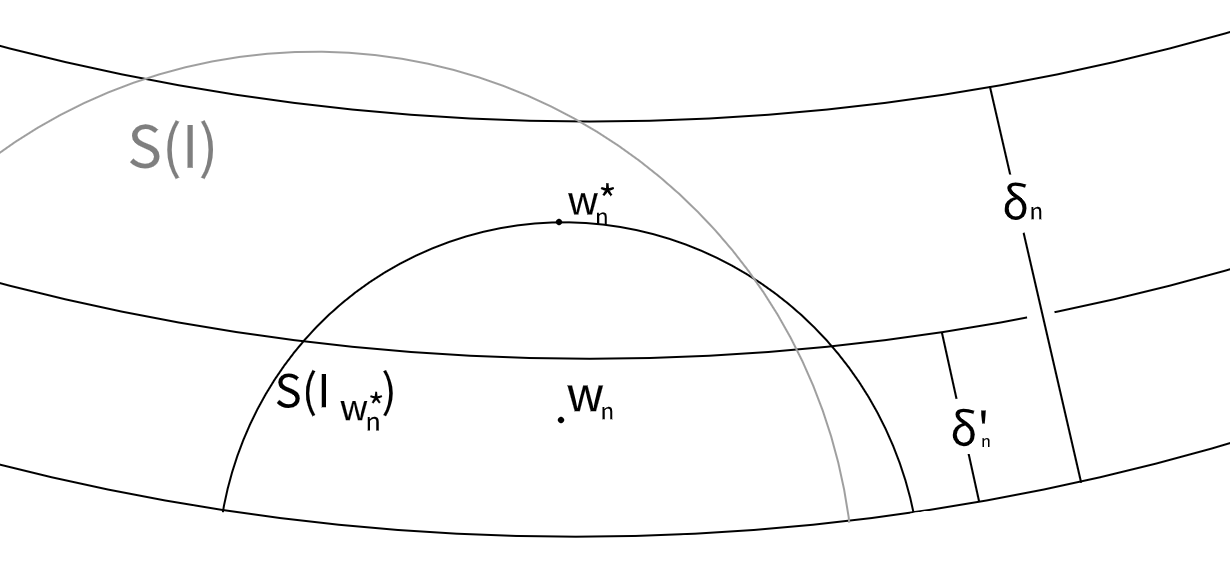}
    \caption{The $n-$th step in the construction }
    \label{fig:recursion}
\end{figure}

\textbf{Recursive definition}
Let $I_0= \partial \bD, w_0=0$ and $a_0=1$. Consequently $F_0\equiv 1$ and conditions $(1)-(3)$ are automatically satisfied. For the recursive step suppose that $I_0,\dots I_{n-1}$, $w_0,\dots w_{n-1}$ and $a_1,\dots a_{n-1}$ are defined and we want to proceed our construction. Since $g\in VMOA$ and $F_{n-1}$ is bounded we can find some $\delta_n >0$ such that 
\begin{equation}\label{eq:smallnearboundary} \sup_{|I|\leq \delta_n} \frac{1}{|I|}\int_{S(I)}|F_{n-1}(z)g'(z)|^2(1-|z|^2)dm(z) \leq 1 . \end{equation}
By part (v) of Lemma \ref{lem:buildingblock} there also exists some $\delta_n'<\delta_n$ such that for any pair of complex numbers $z,w\in \bD$ such that $1-|w| \leq \delta_n'$ and $1-|z| \geq \delta_n$ then $|\beta_w(z)| \leq 3$. We might also choose $\delta'_n$ such that $\sqrt{\delta'_n} \leq 2^{-2n}\delta_n. $

Furthermore, by Lemma \ref{lem:buildingblock}(iv) since $g\not\in BMOA_{\log}$, there exists some $w_n\in \bD $ which satisfies $1-|w_n| \leq \delta_n'$ such that 
\begin{equation}
 \frac{1}{|I_{w_n}|} \int_{I_{w_n}}   \Re (\beta_{w_n}(z))^2|g'(z)|^2(1-|z|^2)dm(z) \geq 2^{2n}.
\end{equation}

Notice furthermore that because $g$ is of vanishing mean oscillation the supremum 
\begin{equation}\label{eq:max} M_n^2:=\sup_{|I| \leq \delta_n}\frac{1}{|I|} \int_{S(I)} \Re (\beta_{w_n}(z))^2|g'(z)|^2(1-|z|^2)dm(z) \geq 2^{2n} \end{equation}
is in fact a maximum which is attained for some interval $I_n \subseteq \partial \bD, |I_n| \leq \delta_n.$ We claim that the function $F_n:=F_{n-1}+M_n^{-1}\beta_{w_n}$ satisfies the required properties.  We start by proving properties $(1)$ and $(2)$. By equation $\eqref{eq:max}$ it is clear that 
\[  a_n:= M_n^{-1} \leq 2^{-n}\]

Also, since $\beta_w$ has positive real part
\[ \frac{1}{|I_n|}\int_{S(I_n)}| F_n(z)g'(z) |^2 (1-|z|^2) dm(z)\geq \frac{M_n^{-2}}{|I_n|}\int_{S(I_n)}|\Re(\beta_{w_n}(z))^2|g'(z)|^2(1-|z|^2)dm(z) = 1.  \]

Finally we need to proved the claimed estimate on $\norm{T_g F_n}_*.$  To obtain this we consider two cases. First suppose that $|I| \geq \delta_n$. We start with a preliminary estimate
\begin{align*}
    \frac{1}{|I|}\int_{S(I)} |\beta_{w_n}(z)g'(z)|^2(1-|z|^2)dm(z) & \leq  \frac{|I_{w_n^*}|}{|I|}\frac{1}{|I_{w_n^*}|}\int_{S(I_{w_n^*})} |\beta_{w_n}(z)g'(z)|^2(1-|z|^2)dm(z) + \\
    & \qquad  \frac{1}{|I|}\int_{S(I)\setminus S(I_{w_n^*})} |\beta_{w_n}(z)g'(z)|^2(1-|z|^2)dm(z) \\ 
    & \lesssim \frac{1-|w_n^*|}{\delta_n} M_n^2 + \norm{g}_*^2 
     \lesssim  \frac{(1-|w_n|)^\frac{1}{2}}{\delta_n}M_n^2+\norm{g}_*^2 \\ & \lesssim \frac{\delta_n'^{\frac{1}{2}}}{\delta_n}M_n^2 +\norm{g}_*^2 \leq 2^{-2n} M_n^2 + \norm{g}_*^2.
\end{align*}
In this estimate we have used the fact that $|I_{w^*}| \approx 1-|w_n^*| \approx (1-|w_n|)^{\frac12}$. Now the induction hypothesis together with the above calculation permit us to estimate as follows
\begin{align*}
 \Big(  \frac{1}{|I|} \int_{S(I)} |F_n(z)g'(z)|^2(1-|z|^2) & dm(z) \Big)^{1/2} \leq \\
 &\leq M_n^{-1} \Big( \frac{1}{|I|}\int_{S(I)}|\beta_{w_n}(z)g'(z)|^2(1-|z|^2) dm(z) \Big)^{1/2} +
 \\ & + \qquad \Big( \frac{1}{|I|}\int_{S(I)} |F_{n-1}(z)g'(z)|^2(1-|z|^2) dm(z) \Big)^{1/2} \\
 &\lesssim 2^{-n} +M_n^{-1} \norm{g}_* +\norm{T_gF_{n-1}}_* \\
 & \leq 2 ^{-n} C(g) + \norm{T_g F_{n-1}}_*,
\end{align*}
where $C(g)$ is a positive constant depending only on $g$ and not on $n$.

It remains to consider the case $ |I| \leq \delta_n$. In this case, equation \eqref{eq:smallnearboundary} allows us to argue as follows;

\begin{align*}
  \Big(  \frac{1}{|I|} &\int_{S(I)} |F_n(z)g'(z)|^2(1-|z|^2) dm(z) \Big)^{1/2}  \leq
  \\
 &\leq\; M_n^{-1} \Big( \frac{1}{|I|}\int_{S(I)}|\beta_{w_n}(z)g'(z)|^2(1-|z|^2) dm(z) \Big)^{1/2} + 
  \\
& \qquad \qquad + \Big( \frac{1}{|I|}\int_{S(I)} |F_{n-1}(z)g'(z)|^2(1-|z|^2) dm(z) \Big)^{1/2}
  \\
&\lesssim \; M_n^{-1} \sup_{|I| \leq \delta_n} \Big(\frac{1}{|I|}\int_{S(I)} \Re(\beta_{w_n}(z))^2|g'(z)|^2(1-|z|^2) dm(z) \Big)^{1/2} + \norm{g}_* + 1
\\
&\; \leq C(g).
\end{align*}

We have proved that 
\[ \norm{T_gF_n}_* \leq \max \{ \norm{T_gF_{n-1}}_*+  2^{-n}C(g), C(g) \}  \]
which completes the induction step and the proof is complete. 
\end{proof}

We also need the Bloch version of Proposition \ref{prop:mainconstruction}

\begin{prop}\label{prop:mainconstructionBloch}
Let $g\in \cB \setminus \cB_{\log, 0}$, then there exists a function $F\in \cB$ such that $T_gF\in \cB \setminus \cB_0 .$
\end{prop}
\begin{proof}
The proof of this proposition is very similar to the proof of Proposition \ref{prop:mainconstruction} and in fact a bit simpler, therefore we shall give only a rough sketch of it. 
A similar argument as the one used in the proof of Proposition \ref{prop:mainconstruction} allows us to reduce the problem to the case $g\in \cB_0 \setminus \cB_{\log}.$

We shall construct inductively two sequences of points in the unit disc $\{z_n\}$ and $\{w_k\}$ such that the functions 
\[ F_n(z) = \sum_{k=0}^na_k \beta_{w_k}(z) \] 
satisfy; 
\begin{itemize}
    \item[(1)] The coefficients $\{ a_k \}$ satisfy $0\leq a_k \leq 2^{-k}.$
    \item[(2)] For all $n\in \bN$ holds
    \[ \Re(F_{n}(z_n))|g'(z_n)|(1-|z_n|^2) \geq 1. \]
    \item[(3)] For all $n \in \bN $ we have $ \norm{T_gF_n}_{\cB} \leq C(g),  $ where $C(g)$ is a constant depending only on $g$.
\end{itemize}

Given this construction we can complete the proof as in the case of functions of bounded mean oscillation. 

Assume without loss of generality that $g'(0) = 1.$ Then set $z_0 = w_0 = 0 $ and $ a_0 = 1 $. For the inductive step, assume that the parameters are defined up to level $n-1$. We can find $\delta_n>0$ such that 
\begin{equation}\label{eq:smallbloch}
\sup_{|z|\geq 1-\delta_m } |F_{n-1}(z)g'(z)|(1-|z|^2) \leq 1. \end{equation}
Choose $\delta_n'>0 $ as in Lemma \ref{lem:buildingblock} (v). Moreover, since $g\not\in \cB_{\log},$ there exists some $w_n\in \bD, 1-|w_n| \leq  \delta_n'$ such that 
\begin{equation}
    \Re\beta_{w_n}(w_n)|g'(w_n)|(1-|w_n|^2) \geq 2^n.
\end{equation}
Finally let $z_n$ a point $1-|z_n|\leq \delta_n$ where the supremum
\[ M_n : = \sup_{1-|z|\leq \delta_n} \Re(\beta_{w_n}(z))|g'(z)|(1-|z|^2) \]
is attained. We finish the recursive step by setting $a_n=M_n^{-1}. $ 
It remains to verify the properties $(1) - (3)$. This is done in a similar way as in the proof of Proposition \ref{prop:mainconstruction} and the details are left to the reader.
\end{proof}

We can now assembly all pieces in order to prove our main result. We shall first prove all equivalences for elliptic semigroups and then we shall prove that no non-elliptic semigroup satisfies $(a)$ or $(b)$ of Theorem \ref{main:thm:elliptic}.

\begin{proof}[Proof of Theorem \ref{main:thm:elliptic}]
We start by proving the equivalence of $(c)$ and $(d)$. Remember that if we assume that $\tau=0$, then
\[
\gamma'(z)\,=\, \frac{z}{G(z)}\,=\,-\frac{1}{p(z)}
\]
but since $\Re(p)\geq 0$, from the Alexander-Noshiro-Warschawski criterion it follows that $\gamma$ is a univalent function in $\cB$, hence in $BMOA$.
Notice that (c) is equivalent to $\gamma \in \cB_{\log, 0}$ and (d) is equivalent to $\gamma \in VMOA_{\log}$, hence the proof is a direct consequence of Corollary \ref{corollary} and the fact that $\gamma$ is univalent.

We proceed now to the proof of the equivalences (a) $\Leftrightarrow$ (c) and (b) $\Leftrightarrow$ (d). The strategy is quite similar for both so we prove in detail that (b) $\Leftrightarrow$ (d) and we sketch the proof for the other implication. 

Let $(\varphi_t)$ be an elliptic semigroup, such that  $[\varphi_t,\,BMOA]=VMOA$. We will show that $\gamma$ must be in $VMOA_{\log}$.
To prove this, suppose that $\gamma \in BMOA\setminus VMOA_{\log}$. From Proposition \ref{prop:mainconstruction}, we can find a function $F\in BMOA$, such that $T_\gamma(F)\in BMOA\setminus VMOA$.
From (\ref{MaxSubTg}) we know that
\[
BMOA \cap \left( T_\gamma(BMOA) \oplus \cC \right) \subseteq [\varphi_t,\,BMOA] \,.
\]
But this means that the function 
$T_\gamma(F) \in  [\varphi_t,\,BMOA] $, which is equal to $VMOA$ by our assumption, and this is a contradiction.

Conversly, suppose that (d) holds, i.e. the infinitesimal generator $G$ satisfies  (\ref{eq:lvmo}),  which also implies (\ref{eq:lbmo}). Then from \cite[Corollary 2]{Blasco2013} the result follows.

For the equivalence of (a) and (c), one needs to follow the exact same reasoning as before. In specific, assuming that (a) holds, use Proposition \ref{prop:mainconstructionBloch} together with the fact that
\[
\cB \cap \left( T_\gamma(\cB) \oplus \cC \right) \subseteq [\varphi_t,\,\cB] \,
\]
to ensure that $\gamma \in \cB_{\log,0}$ by contradiction, and for the converse implication apply \cite[Corollary 2]{Blasco2013} as before.

It remains to prove that non-elliptic semigroups  cannot satisfy $[\varphi_t,BMOA]= VMOA$ or $[\varphi_t, \cB]=\cB_0$.
The argument is identical for both $BMOA$ and $\cB$ so we shall only deal with the Bloch space. 
Assume that $(\varphi_t)$ is a non-elliptic semigroup such that $[\varphi_t,\cB]=\cB_0$. Without loss of generality we consider $\tau = 1 $. Let also $h$ be the associated Koenigs function and $H(z)=h(z)/z$. Then classical Koebe's distortion theorem implies that  $\log H \in \cB$ (in fact by \cite[p. 87, Remark 2]{Blasco2013} we know that  $\log H \in\cB_0$ but we shall not need this extra information). Then we distinguish two cases. If $\log H \in \cB_{\log,0}$ then we know from Theorem \ref{bound.Comp.} that the operator 
\[ T_{\log H}: \cB \to \cB \]
is {\it compact}. Then if $\lambda\neq 0$ is a point in the spectrum of $T_{\log H}$, by the spectral theorem for compact operators \cite[Section 21.2]{lax2002functional} must be an eigenvalue. But this is impossible since $T_{\log H} f = \lambda f$ implies that $f\equiv 0$ \cite[Proposition 5.1]{Aleman2009}. Therefore $T_{\log H}$ has trivial spectrum. In particular there exists $f \in \cB$ such that 
\[f(z) - T_{\log H} f (z)  \equiv 1. \]
Solving this first order ODE we find that $f = H$, which implies that $h \in \cB  $, or equivalently by Pommerenke's Theorem \cite{Pommerenke1978} $h \in BMOA$. This contradicts \eqref{not BMOA}.

This leaves only the possibility that $\log H \in \cB \setminus \cB_{\log,0}.$ Then we are again in a situation where we can apply Proposition \ref{prop:mainconstructionBloch}. Therefore there exists $F\in \cB$ such that $T_{\log H} F  \in \cB \setminus \cB_0. $ Notice that this is equivalent to the fact that the function
\[ T_h\Big(\frac{F}{H}\Big)(z) = \int_0^z \frac{tF(t)h'(t)}{h(t)} dt    \]
 belongs to $\cB \setminus \cB_0$. Since in this case $h  $ is the $\gamma$-symbol of the semigroup it remains to prove that the function $F/H$ is a Bloch function in order to arrive at a contradiction.
 
 We have,
 \begin{align*}
     \Big| \Big( \frac{F}{H} \Big)'(z) \Big| (1-|z|^2) & \leq \frac{|F'(z)|(1-|z|^2)}{|H(z)|} + \frac{|F(z)H'(z)|(1-|z|^2)}{|H(z)|^2} \\ 
     & \leq \norm{H^{-1}}_{H^{\infty}} \big( \norm{F}_{\cB} + \norm{T_{\log H} F}_{\cB} \big)<\infty.
 \end{align*}
 In other words, we have shown that the function $ T_h\Big(\frac{F}{H}\Big) \in \cB$ and at the same time it belongs to the range $T_\gamma(\cB)$, hence by \eqref{MaxSubTg} in $[\varphi_t, \cB]=\cB_0$ and this is a contradiction.
\end{proof}

\subsection*{Further remarks} We believe that the techniques that we have employed can be used to prove similar characterizations in other kind of spaces. In particular, recent studies investigated the maximal subspace of holomorphic semigroups in $BMOA$-type spaces \cite{Daskalogiannis2021}, in the analytic Morrey spaces \cite{Galanopoulos2020} and in the setting of mixed norm spaces \cite{Arevalo2019}. It would be interesting to know whether a similar characterization of the minimality of the maximal subspace is possible in these settings, too.

\bibliography{literature.bib}
\bibliographystyle{plain}

\end{document}